\documentclass[a4paper,11pt]{article}
\usepackage{latexsym}
\usepackage{amssymb}
\usepackage{amsfonts}
\usepackage{amsmath}
\usepackage{indentfirst}
\newtheorem{prop}{Proposition}[section]
\newtheorem{teor}{Theorem}[section]

\newtheorem{cor}{Corollary}[section]
\newcommand{\nkinN}{n,k\in \mathbf{N}}
\newcommand{\ninN}{n\in \mathbf{N}}
\newcommand{\nkinNstar}{n,k\in \mathbf{N^{*}}}
\newcommand{\ninNstar}{n\in \mathbf{N^{*}}}
\newcommand{\kinNstar}{k\in \mathbf{N^{*}}}
\newcommand{\cvd}{\quad $\blacksquare$\bigskip}
\newcommand{\Rx}{\mathbf{R}[x]}
\date{}

\author{Luca Ferrari\thanks{Dipartimento di Scienze Matematiche ed
Informatiche, Pian dei Mantellini, 44, 53100, Siena, Italy {\tt
ferrari@math.unifi.it}} \and Renzo Pinzani\thanks{Dipartimento di
Sistemi e Informatica, via Lombroso 6/17, 50134 Firenze, Italy
{\tt pinzani@dsi.unifi.it}}}

\title{Catalan-like numbers and succession rules\footnote{This work was partially supported
by MIUR project: \emph{Linguaggi formali e automi: metodi, modelli
e applicazioni}.}}

\begin{document}

\maketitle

\begin{abstract}
The ECO method and the theory of Catalan-like numbers introduced
by Aigner seems two completely unrelated combinatorial settings.
In this work we try to establish a bridge between them, aiming at
starting a (hopefully) fruitful study on their interactions. We
show that, in a linear algebra context (more precisely, using
infinite matrices), a succession rule can be translated into a
(generalized) Aigner matrix by means of a suitable change of basis
in the vector space of one-variable polynomials. We provide some
examples to illustrate this fact and apply it to the study of two
particular classes of succession rules.
\end{abstract}

\begin{center}
\end{center}

\section{Introduction}

The ECO method was founded in the 90's by a group of researchers,
including Pinzani, Barcucci, Del Lungo and Pergola
\cite{BDLPP,BDLPP1}. It consists of a purely combinatorial way of
constructing the objects of a given class in such a way that, if
the construction is sufficiently regular and recursive,
enumeration follows by more or less standard methods of
combinatorial analysis. More precisely, one starts by partitioning
a class of objects according to their size (suitably defined). The
goal is then to perform a sort of local expansion on each object
of a given size, thus producing all the objects of the successive
size exactly once. Therefore a single object produces a set of new
objects according to some parameter. Typically, if such a
construction is regular enough, one can encode it using a
\emph{succession rule} \cite{W1,W2}, which is a purely formal
system, generally expressed as follows:
\begin{equation}\label{regola}
\left\{ \begin{array}{ll} (a)
\\ (k)\rightsquigarrow (e_1 (k))\cdots (e_k (k))
\end{array}
\right. .
\end{equation}

Here the letters denote positive integers, $(a)$ is called the
\emph{axiom} and $(k)\rightsquigarrow (e_1 (k))\cdots (e_k (k))$
is the \emph{production} of $(k)$. A succession rule can be
represented by means of its \emph{generating tree}, which is, by
definition, the infinite, rooted, labelled tree whose root is
labelled $a$ (like the axiom) and such that every node labelled
$k$ produces $k$ sons, labelled respectively $e_1 (k),\ldots ,e_k
(k)$. One of the main enumerative information provided by a
succession rule is the numerical sequence of the cardinalities of
the levels of the generating tree associated with the rule: we
will refer to such a sequence as the numerical sequence
\emph{determined by} the succession rule.

The basic reference for the ECO method is \cite{BDLPP}, in which
many examples can also be found.

The importance of succession rules as a tool for the ECO method
has lead to several investigations to get a better mathematical
insight on them. In \cite{FP} the authors define the concept of
\emph{rule operator}, thus translating a succession rule into a
linear operator on one-variable polynomials. In \cite{DFR} every
succession rule is associated with at least two infinite matrices:
the \emph{production matrix}, which is essentially the matrix of
the related rule operator with respect to the canonical basis of
polynomials $(x^n )_{\ninN}$, and the \emph{ECO matrix}, whose
$(n,k)$-entry is, by definition, the number of nodes labelled $k$
at level $n$ in the corresponding generating tree. We point out
that a few instances of the notion of production matrix appeared
for the first time in \cite{W1} under the name of ``transfer
matrices".

\bigskip

Another combinatorial theory dealing with infinite matrices is
Aigner's theory of Catalan-like numbers \cite{A1,A2,A3,Z}. The
basic idea of \cite{A1} is to characterize those sequences for
which the determinants of the Hankel matrices of order 0 are all
equal to 1. It is shown that such sequences appear as the first
column of certain infinite matrices, called \emph{admissible
matrices}. These numbers are referred to with the name of
\emph{Catalan-like numbers}. The reason for this name lies in the
fact that Catalan numbers are the unique sequence whose Hankel
determinants of orders 0 and 1 equals 1.

In \cite{A2,A3} Aigner extends this theory by considering a more
general kind of matrices, which we will rename \emph{Aigner
matrices} (instead of the infelicitous name ``recursive matrices"
used in \cite{A3}). Generalizing the previous definition, we will
call \emph{Catalan-like numbers} every sequence appearing as the
first column of an Aigner matrix.

\bigskip

The aim of our work is to provide a ``vocabulary" to translate the
ECO method into Aigner's theory, and vice versa. Such a vocabulary
turns out to be based on linear algebra tools, consisting of a
suitable change of basis in the vector space of one-variable
polynomials. What we hope to show in this paper is that the two
theories under consideration are, in some algebraic sense, the two
sides of the same medal, which is quite surprising if we think of
the very different starting points, and combinatorial meanings, of
such theories.

After a brief survey of the notions we need from the ECO method
and Aigner's theory, we provide the main linear algebra tools to
be used in the sequel. In particular, we define the Aigner basis
in the vector space of one-variable polynomials and prove some of
its properties. Next we introduce what we call the
\emph{fundamental change of basis}, which is the key ingredient to
accomplish our project, and propose some examples to illustrate
our approach. The final part of the work is devoted to the study
of two particular cases, for which we are able to fully describe
how to switch from one theory to the other. In the last section we
give some hints for possible, future works.

\bigskip

Before starting, a few words concerning our notations. We have
chosen to use $\mathbf{N}$ to denote natural numbers (zero
included), whereas $\mathbf{N}^{*}$ is the set of natural numbers
without zero. The symbol $\mathbf{x}$ is used to denote the
operator of multiplication by $x$, in order to distinguish it from
the symbol $x$, used as a polynomial variable. The symbol $^\top$
denotes the transpose of a matrix. The last remark concerns the
way we have chosen to index the lines of our matrices.
Classically, the lines of an Aigner matrix are indexed by
$\mathbf{N}$, whereas, in an ECO matrix, the first column is
usually column 1 (so that columns are indexed by $\mathbf{N}^*$).
There are clear combinatorial and algebraic reasons for this: the
$n$-th column of an ECO matrix gives the distribution of label
$(n)$ in the generating tree, whereas the scalar product of the
$n$-th and the $m$-th rows of an admissible matrix gives the first
element of its $(n+m)$-th row. Unfortunately, keeping both these
conventions would result in a purely formal, but greatly
inelegant, variation of our results: namely, the fundamental
change of basis (which is degree-preserving in our theory) would
translate $x^n$ into a polynomial of degree $n-1$ (which is
$p_{n-1}(x)$, according to the notations of section \ref{tools}).
To avoid this we have preferred to use $\mathbf{N}^{*}$ as the set
of indices for the lines of Aigner matrices. In this way the
(nice) characteristic algebraic property of admissible matrices
becomes a little bit difficult to read, but our theory can be
described much more elegantly.

\section{Preliminaries on ECO and Aigner
matrices}\label{preliminari}

In this section we report the main facts concerning the two
combinatorial theories we are going to compare.

\bigskip

Consider a succession rule as in (\ref{regola}). Instead of
representing it by means of a generating tree, one can choose
linear algebra notations. In the vector space of one-variable
polynomials, define the linear operator $L=L_{\Omega}$ on the
canonical basis $(x^n )_{\ninN}$ as follows:
\begin{eqnarray*}
L(\mathbf{1})&=&x^a ,
\\ L(x^k )&=&x^{e_1 (k)}+\ldots +x^{e_k (k)},
\\ L(x^h )&=&hx^h ,\qquad \textnormal{if $(h)$ is not a label of $\Omega$.}
\end{eqnarray*}

$L$ is called the \emph{rule operator} associated with $\Omega$
(see \cite{FP}): it bears all the enumerative properties of a
succession rule and allows to express such properties using
algebraic notations. For example, if $(f_n )_{\ninN}$ is the
numerical sequence determined by $\Omega$, then we can find $f_n$
using the rule operator $L$ of $\Omega$ as follows:
\begin{displaymath}
f_n =[L^{n+1}(\mathbf{1})]_{x=1}.
\end{displaymath}

Throughout the present work, we will always deal with a special
case, namely we assume that deg $L^n (\mathbf{1})=n$. From a
combinatorial point of view, this means that the set of the labels
of $\Omega$ is $\mathbf{N}^{*}$ and the maximum label among those
produced by $(k)$ is $(k+1)$.

The infinite matrix $P=P_{\Omega}$ representing $L$ with respect
to the canonical basis $(x^n )_{\ninN}$ is called the
\emph{production matrix} of $\Omega$. Such matrices are
extensively studied in \cite{DFR}; here we only recall some of
their properties.

Let $A_P$ be the infinite matrix whose $n$-th row vector is
given by $u^\top P^{n-1}$ (where $u^\top =(1\; 0\; 0\; \ldots \;
0\; \ldots )$). Then $A_P$ describes the statistic given by the
distribution of the labels at the various levels of the generating
tree related to $\Omega$. Namely, the $(n,k)$ entry of $A_P$ is
the number of nodes labelled $k$ at level $n$ of the generating
tree of $\Omega$. $A_P$ is called the \emph{ECO matrix} associated
with $P$ (or with $\Omega$), and is also characterized by the
matrix equality $DA_P =A_P P$, where
\begin{displaymath}
D=\begin{pmatrix} 0 & 1 & 0 & 0 & 0 & \ldots
\\ 0 & 0 & 1 & 0 & 0 & \ldots
\\ 0 & 0 & 0 & 1 & 0 & \ldots
\\ 0 & 0 & 0 & 0 & 1 & \ldots
\\ 0 & 0 & 0 & 0 & 0 & \ldots
\\ \vdots & \vdots & \vdots & \vdots & \vdots & \ddots
\end{pmatrix}.
\end{displaymath}
In terms of the production matrix $P$, the sequence $(f_n
)_{\ninN}$ is nothing else than the sequence of the row sums of
the associated ECO matrix.
The ordinary and exponential
generating functions of $\Omega$ are given, respectively, by $f_P
(t)=u^\top (I-tP)^{-1}e$ and $F_P (t)=u^\top \exp{(tP)}e$, where
$e$ is the column vector $(1\; 1\; 1\; \ldots 1\; \ldots )^\top$
and $\exp{(X)}$ denotes the usual matrix exponential of the
(infinite) matrix $X$.
In section \ref{esempi} we deal with some examples of the theory
we are going to develop; in describing such examples we also
consider the production and ECO matrices of some classical
succession rules. The reader is invited to have a look to those
examples in order to be introduced to these concepts.

\bigskip

In a recent series of nice and well-written papers
\cite{A1,A2,A3}, Martin Aigner has developed a new theory to deal
with numerical sequences somehow linked to the sequence of Catalan
numbers. One of the main tool of this theory is a particular class
of infinite (triangular) matrices, called \emph{admissible
matrices} in \cite{A1} and renamed with the infelicitous term
\emph{recursive matrices} after \cite{A3}. In the sequel we will
use the following terminology.

\bigskip

Consider a lower triangular matrix $A=(a_{nk})_{\nkinNstar}$ with
main diagonal equal to 1. $A$ is called an \emph{admissible
matrix} whenever, for every $n,m\in \mathbf{N^*}$, the (ordinary)
scalar product of the $n$-th and the $m$-th rows of $A$ gives the
first element of the $(n+m-1)$-st row; in symbols:
\begin{displaymath}
\sum_{k\geq 1} a_{nk}a_{mk}=a_{(n+m-1)1}.
\end{displaymath}

More generally, we define $A$ to be an \emph{Aigner matrix} when
there exists a sequence of nonnegative integers $(T_n
)_{\ninNstar}$ such that $T_1 =1$ and $T_n |T_{n+1}$ for which we
have
\begin{equation}\label{caraig}
\sum_{k\geq 1} a_{nk}a_{mk}T_k =a_{(n+m-1)1}.
\end{equation}

Obviously admissible matrices correspond to Aigner matrices for
which $T_n \equiv 1$.

If $A$ is an Aigner matrix, the sequence $(a_{n1})_{\ninNstar}$ is
called the sequence of \emph{Catalan-like numbers} associated with
$A$. For an extensive study of the algebraic and enumerative
properties of Aigner matrices the reader is referred to
\cite{A1,A3} and to the further items cited in the references of
the two papers. Here we recall only those results which we need
for our purposes.

\begin{prop} (\cite{A1,A3}) An Aigner matrix $A=(a_{nk})_{\nkinNstar}$ is uniquely
determined by the two sequences $(a_{(n+1)n})_{\ninNstar}$ and
$(T_n )_{\ninNstar}$. Conversely, to every pair of sequences $(b_n
)_{\ninNstar}$ and $(T_n )_{\ninNstar}$ of real numbers there
exists an (and therefore precisely one) Aigner matrix
$A=(a_{nk})_{\nkinNstar}$ associated with $(T_n )_{\ninNstar}$ and
such that $a_{(n+1)n}=b_n$ for all $n$.
\end{prop}

\begin{prop}\label{boh} (\cite{A1,A3}) Let $A=(a_{nk})_{\nkinNstar}$ be an Aigner matrix
associated with $(T_n )_{\ninNstar}$. Set $s_1 =a_{21},s_n
=a_{(n+1)n}-a_{n(n-1)}$ and $t_n =\frac{T_n}{T_{n-1}}$, for $n\geq
2$. Then we have
\begin{eqnarray}\label{aigrec}
a_{11}&=&1,\qquad \qquad a_{1k}=0 \quad (k>1)\nonumber
\\ a_{nk}&=&a_{(n-1)(k-1)}+s_{k}a_{(n-1)k}+t_{k+1}a_{(n-1)(k+1)}\quad (n\geq
2).
\end{eqnarray}

Conversely, if $a_{nk}$ is given by the recursion (\ref{aigrec}),
then $A=(a_{nk})_{\nkinNstar}$ is an Aigner matrix with $T_n =t_2
\cdot \ldots \cdot t_n$ and $a_{(n+1)n}=s_1 +\ldots +s_n$.
\end{prop}

Setting $\sigma =(s_n )_{\ninNstar}$, $\tau =(t_n )_{n\geq 2}$, we
say that $A=A^{\sigma ,\tau}$ is the Aigner matrix \emph{of type}
$(\sigma ,\tau )$ when recursion (\ref{aigrec}) holds for its
entries.

It is possible to write (\ref{caraig}) in a compact matrix form.
Define the diagonal matrix $T$ and the infinite Hankel matrix of
the sequence $(a_{n1})_{\ninNstar}$ as follows:
\begin{displaymath}
T=\left( \begin{array}{llll} T_1&0&0&\cdots
\\ 0&T_2&0&\cdots
\\ 0&0&T_3&\cdots
\\ \vdots &\vdots &\vdots &\ddots
\end{array}\right) ,\quad
H=\left( \begin{array}{llll} a_{11}&a_{21}&a_{31}&\cdots
\\ a_{21}&a_{31}&a_{41}&\cdots
\\ a_{31}&a_{41}&a_{51}&\cdots
\\ \vdots &\vdots &\vdots &\ddots
\end{array}\right) .
\end{displaymath}

Then (\ref{caraig}) can be written as $ATA^\top =H$. More
precisely, we have the following characterization.

\begin{prop} (\cite{A3}) $A$ is an Aigner matrix if and only if $ATA^\top =H$
with $T_n \neq 0$ for all $n\geq 2$, $T_1 =1$. The sequences
$\sigma$ and $\tau$ are then given as in proposition \ref{boh}.
\end{prop}

Moreover, if we denote by $H_n$ the $n$-th Hankel matrix of a
sequence $(a_{n})_{\ninNstar}$ (which is, by definition, the
submatrix of $H$ consisting of rows and columns 1 to $n$), we have
the following corollary.

\begin{cor} (\cite{A3}) A sequence $(a_{n})_{\ninNstar}$ is
Catalan-like if and only if $|H_n|\neq 0$ for all $n\geq 1$.
\end{cor}

\section{Linear algebra tools}\label{tools}

In this section we give some (elementary) linear algebra tools
which will be necessary in the sequel to settle the stated analogy
between ECO method and Aigner's theory of Catalan-like numbers.

\bigskip

In the vector space of one-variable polynomials over the real
field (to be denoted $\Rx$) we define the following polynomial
sequence:
\begin{eqnarray*}
p_0 (x)&=&1,
\\ p_n (x)&=&x(x-1)^{n-1},\qquad \forall n\geq 1.
\end{eqnarray*}

It is clear that deg $p_n (x)=n$, so that $(p_n (x))_{\ninN}$
constitutes a basis for the vector space $\Rx$. We call such a
basis the \emph{Aigner basis} of $\Rx$.

\bigskip

\emph{Remark.} We recall that the polynomial $p_n (x)$ has a very
important combinatorial meaning: it is the chromatic polynomial of
a tree having $n$ vertices (see, for example, \cite{B}). However,
this fact will not be used in this paper.

\bigskip

It is well known (see, for instance, \cite{A0}) that, for any
basis of $\Rx$, there exists a unique differential\footnote{i.e.,
it maps a polynomial of degree $n$ into a polynomial of degree
$n-1$.} linear operator mapping the $n$-th element of the basis
into the $(n-1)$-st one. According to \cite{FP}, we call
\emph{factorial derivative operator} the linear operator:
\begin{eqnarray*}
T&:&\Rx\longrightarrow \Rx
\\ &:&p_0 (x)\longrightarrow \mathbf{0},
\\ &:&p_n (x)\longrightarrow p_{n-1}(x),\qquad n\geq 1.
\end{eqnarray*}

The factorial derivative operator can also be nicely performed on
the canonical basis of $\Rx$.

\begin{prop} For any $\ninN$, $T(x^n )=1+x+\ldots
+x^{n-1}=\sum_{k=0}^{n-1}x^k$.
\end{prop}

\emph{Proof.} For the first values of $n$, we have
$T(\mathbf{1})=\mathbf{0}$, $T(x)=\mathbf{1}$, $T(x^2 )=T(p_1
(x)+p_2 (x))=1+x$. By induction, suppose that $T(x^n )=1+\ldots
x^{n-1}$. Observe that we have the trivial equality $T((x-1)p_n
(x))=T(p_{n+1}(x))=p_n (x)$, whence $T(xp_n (x))=T(p_n (x))+p_n
(x)$. By linearity we then have:
\begin{displaymath}
T(xp(x))=T(p(x))+p(x).
\end{displaymath}

Therefore, in the case $p(x)=x^n$, we get:
\begin{eqnarray*}
T(x^{n+1})&=&T(x\cdot x^n )=T(x^n )+x^n
\\ &=&1+\ldots +x^{n-1}+x^n ,
\end{eqnarray*}
which is the thesis.\cvd

\bigskip

Some properties of the factorial derivative operator are recorded
in \cite{FP}.

%
%
%
\bigskip

The Aigner basis has a nice behavior with respect to the usual
multiplication operation on polynomials.

\begin{prop} Let $n,m\geq 1$.
\begin{enumerate}
\item $p_n (x)\cdot p_m (x)=xp_{n+m-1}(x)$; \item $x^k p_n
(x)=\sum_{h=0}^{k}{k\choose h}p_{n+h}(x)$; in particular, $xp_n
(x)=p_{n+1}(x)+p_n (x)$; \item denoting by $\mathbf{x}^{-1}$ the
linear operator defined by
$\mathbf{x}^{-1}(p(x))=\frac{p(x)-p(0)}{x}$ (so that
$\mathbf{x}^{-1}$ is the usual difference quotient operator), it
is $\mathbf{x}^{-k}p_n
(x)=\sum_{h=0}^{n-k}(-1)^{n-k-h}{n-1-h\choose k-1}p_h (x)$; in
particular, setting $k=1$, we have $\mathbf{x}^{-1}p_n
(x)=\sum_{h=0}^{n-1}(-1)^{n-h-1}p_h
(x)=p_{n-1}(x)-p_{n-2}(x)+p_{n-3}(x)-\ldots$.
\end{enumerate}
\end{prop}

\emph{Proof.}
\begin{enumerate}
\item $p_n (x)\cdot p_m (x)=x(x-1)^{n-1}\cdot x(x-1)^{m-1}=x\cdot
x(x-1)^{n+m-2}=xp_{n+m-1}(x)$. \item We have immediately
\begin{displaymath}
p_{n+1}(x)+p_n (x)=x(x-1)^{n-1}(x-1+1)=xp_n (x).
\end{displaymath}

Then, by induction, we get
\begin{eqnarray}\label{formula}
x^{k+1}p_n (x)&=&x\cdot x^k p_n (x)=x\sum_{h=0}^{k}{k\choose
h}p_{n+h}(x)\nonumber
\\ &=&\sum_{h=0}^{k}{k\choose h}(p_{n+h+1}(x)+p_{n+h}(x))\nonumber
\\ &=&p_n (x)+\sum_{h=1}^{k}\left( {k\choose h}+{k\choose
h-1}\right) p_{n+h}(x)+p_{n+k+1}(x)\nonumber
\\ &=&\sum_{h=0}^{k+1}{k+1\choose h}p_{n+h}(x).
\end{eqnarray}
\item Just recall that $\mathbf{x}^{-k}$ is the inverse of the
operator $\mathbf{x}^k$. Then the thesis is obtained by simply
inverting the combinatorial sum in (\ref{formula}).\cvd
\end{enumerate}

We close this section by stating a technical result, useful in the
computation of the powers of the factorial derivative operator
$T$, which can be proved by induction.

\begin{prop} For $n,k\in \mathbf{N}$, we have:
\begin{displaymath}
T^k (x^n )=\sum_{h=0}^{n-k}{n-h-1\choose k-1}x^h .
\end{displaymath}
\end{prop}

\begin{cor}\label{powers} \begin{enumerate}
\item $[T^k (x^n )]_{x=1}={n\choose k}$; \item $[T^k (x^n
)]_{x=0}={n-1\choose k-1}$.
\end{enumerate}
\end{cor}

\section{The fundamental change of basis}\label{bascha}

For our purposes, we slightly modify the definition of ECO matrix
given in \cite{DFR}, namely we suppose that the $n$-th row of $F$
(with $\ninNstar$) gives the distribution of the various labels
\emph{at level $n-1$} of the generating tree of the rule; thus the
$(n,k)$ entry of $F$ is the number of nodes labelled $k$ \emph{at
level $n-1$}. Let $\Omega$ be a succession rule as in
(\ref{regola}) and suppose that $F=(f_{nk})_{\nkinNstar}$ is the
ECO matrix associated with $\Omega$, as explained above. We denote
by $\beta_n (x)$ the polynomial canonically associated with the
$n$-th row of $F$, namely:
\begin{equation}\label{rowpoly}
\beta_n (x)=\sum_{k=1}^{n}f_{nk}x^k \qquad (n\geq 1).
\end{equation}

As we have already remarked, we will always assume that deg
$\beta_n (x)=n$.

Now expand the polynomials $\beta_n (x)$ in terms of the Aigner
basis, thus obtaining
\begin{equation}\label{change}
\beta_n (x)=\sum_{k=1}^{n} a_{nk}p_k (x).
\end{equation}

Clearly, the coefficients $f_{nk}$ and $a_{nk}$ are intimately
related. In particular, the following, simple result shows that
this setting is the right one to achieve our project.

\begin{prop} The succession $(r_n)_{\ninNstar}=(
\sum_{k=1}^{n}f_{nk})_{\ninNstar}$ of the row sums of $F$ is equal
to the succession $(a_{n1})_{\ninNstar}$ given by the first column
of the matrix $A=(a_{nk})_{\nkinNstar}$. In symbols:
\begin{displaymath}
a_{n1}=\sum_{k=1}^{n}f_{nk}.
\end{displaymath}
\end{prop}

\emph{Proof.} It is obvious that $r_n =\beta_n (1)$, whence
$\sum_{k=1}^{n}a_{nk}p_k (1)=r_n$. By the definition of the Aigner
basis, it is $p_k (1)=0$, for $k>1$, and $p_1 (1)=1$, and so $r_n
=a_{n1}$, as desired.\cvd

\bigskip

The results obtained so far can be naturally expressed also in
matrix notation. Specifically, it turns out that the fundamental
change of basis described in (\ref{change}) is represented as the
multiplication on the right by the Pascal matrix. In other words,
if $P=\left( {n\choose k}\right)_{\nkinN}$ is the usual Pascal
matrix, and $F=(f_{nk})_{\nkinNstar}$ a given ECO matrix, then the
associated matrix $A=(a_{nk})_{\nkinNstar}$ can be expressed as
follows:
\begin{displaymath}
A=FP.
\end{displaymath}

So the Pascal matrix $P$ is the matrix of the change of basis from
$(x^n )_{\ninN}$ to the Aigner basis $(p_n (x))_{\ninN}$.

\bigskip

Now cards are laid on the table: changing the canonical basis into
the Aigner basis is the ``linear algebra" way to switch from the
ECO method to Aigner's theory. At this stage, the following, very
natural question can be asked:

\bigskip

\emph{1) for which ECO matrices $F$ does it happen that the matrix
$A$ is an Aigner matrix?}

\bigskip

Regarding things the other way round, one can start with an Aigner
matrix $A$ and perform the inverse change of basis (from $p_n (x)$
to $x^n$). Therefore the previous question can be inverted:

\bigskip

\emph{2) for which Aigner matrices $A$ does it happen that the
matrix $F$ is an ECO matrix?}

\bigskip

These two problems seems to be rather difficult to be tackled in
their full generality. In the rest of the paper we will mainly
focus on special examples to hopefully illustrate the interest of
our approach. Finally, we will consider two particular classes of
ECO matrices (namely, those arising from the so-called
\emph{factorial succession rules} and \emph{differential
succession rules}), giving for them a complete answer to question
1.

\bigskip

Before closing this section, we record some further notations and
results which will be useful in the sequel.

If $L$ is the rule operator associated with $\Omega$, then for the
polynomials $\beta_n (x)$ in (\ref{rowpoly}) we clearly have
\begin{displaymath}
\beta_n (x)=L^n (\mathbf{1}).
\end{displaymath}

Applying $L$ means to shift from row $n$ to row $n+1$ in the ECO
matrix associated with $\Omega$, whence:
\begin{equation}\label{operiga}
L(\beta_n (x))=\beta_{n+1}(x).
\end{equation}

The coefficients $f_{nk}$ and $a_{nk}$ in (\ref{rowpoly}) and
(\ref{change}) can be expressed in linear algebraic terms, as the
following proposition clarifies.

\begin{prop}\label{coefficienti}
\begin{enumerate}
\item $f_{nk}=\left[ \frac{D^k}{k!}(\beta_n (x))\right]_{x=0}$.
\item $a_{nk}=\left[ T^k (\beta_n (x))\right]_{x=0}$.
\end{enumerate}
\end{prop}

The easy proof is left to the reader.

\begin{cor}\label{hint} For any $n$,
$a_{(n+1)n}=f_{(n+1)n}+nf_{(n+1)(n+1)}$.
\end{cor}

\emph{Proof.} From the above proposition we have
\begin{displaymath}
a_{(n+1)n}=[T^n (\beta_{n+1}(x))]_{x=0}=\sum_{k\geq
1}f_{(n+1)k}[T^n (x^k )]_{x=0}.
\end{displaymath}

Now, recalling corollary \ref{powers}, we get immediately:
\begin{displaymath}
a_{(n+1)n}=\sum_{k\geq 1}{k-1\choose
n-1}f_{(n+1)k}=f_{(n+1)n}+nf_{(n+1)(n+1)},
\end{displaymath}
as desired.\cvd

To conclude, we prove the following proposition, concerning the
behavior of a rule operator when applied to the Aigner basis.

\begin{prop} For any rule operator $L$ and for any $n>2$, we have:
\begin{displaymath}
[L(p_n (x))]_{x=1}=0.
\end{displaymath}
\end{prop}

\emph{Proof.} It is easy to see \cite{FP} that, for any rule
operator $L$, it is
\begin{displaymath}
[L(p(x))]_{x=1}=[D(p(x))]_{x=1}.
\end{displaymath}

(This is due to the fact that $[L(x^k )]_{x=1}=k=[D(x^k
)]_{x=1}$). Then we get immediately:
\begin{displaymath}
[L(p_n (x))]_{x=1}=[D(p_n
(x))]_{x=1}=[Dx(x-1)^{n-1}]_{x=1}=0.\quad \blacksquare
\end{displaymath}

\section{Detailed examples}\label{esempi}

In this section we will provide a detailed analysis of the effects
of the fundamental change of basis in the case of a well-known
succession rule determining Catalan numbers. Then some other
examples will be dealt with; for them, we will only state the main
facts (without giving proofs), however our results can be checked
out by a direct computation or by applying the theory we are going
to develop in the final part of our work.

\subsection{Catalan numbers}

Let us consider the succession rule
\begin{displaymath}
\Omega: \left\{ \begin{array}{lll} (1)
\\ (1)\rightsquigarrow (2)
\\ (k)\rightsquigarrow (2)(3)(4)\cdots (k)(k+1)
\end{array}
\right. ,
\end{displaymath}
defining Catalan numbers 1,1,2,5,14,42,132,$\ldots$ (see, for
example, \cite{BDLPP}). The first lines of the ECO matrix
associated with $\Omega$ looks as follows:
\begin{displaymath}
F=\left( \begin{array}{lllllllll} 1&0&0&0&0&0&0&0&\cdots
\\ 0&1&0&0&0&0&0&0&\cdots
\\ 0&1&1&0&0&0&0&0&\cdots
\\ 0&2&2&1&0&0&0&0&\cdots
\\ 0&5&5&3&1&0&0&0&\cdots
\\ 0&14&14&9&4&1&0&0&\cdots
\\ 0&42&42&28&14&5&1&0&\cdots
\\ 0&132&132&90&48&20&6&1&\cdots
\\ \vdots &\vdots &\vdots &\vdots &\vdots &\vdots &\vdots &\vdots
&\ddots
\end{array}\right) .
\end{displaymath}

Applying the fundamental change of basis, one gets the following
matrix:
\begin{displaymath}
A=\left( \begin{array}{llllll} 1&0&0&0&0&\cdots
\\ 1&1&0&0&0&\cdots
\\ 2&3&1&0&0&\cdots
\\ 5&9&5&1&0&\cdots
\\ 14&28&20&7&1&\cdots
\\ \vdots &\vdots &\vdots &\vdots &\vdots &\ddots
\end{array}\right) .
\end{displaymath}

Is $A$ the (unique) admissible matrix for Catalan numbers found in
\cite{A1}? The first fact suggesting that it could be so follows
from the application of corollary \ref{hint}. It is well known
\cite{FP} that the entries of $F$ are the so-called \emph{ballot
numbers}, namely:
\begin{displaymath}
f_{nk}=\frac{k-1}{n-1}{2n-k-2\choose n-k}\qquad (n,k>1).
\end{displaymath}

Therefore it follows immediately that, in $A$, we have:
\begin{eqnarray*}
a_{(n+1)n}&=&f_{(n+1)n}+nf_{(n+1)(n+1)}
\\ &=&\frac{n-1}{n}{n\choose 1}+n\cdot \frac{n}{n}{n-1\choose
0}=n-1+n=2n-1,
\end{eqnarray*}
which agrees with the $(n+1,n)$ entry of the admissible matrix for
Catalan numbers. Obviously, this is not enough to conclude,
however it is in fact a strong hint. To get to the desired result
we need to show, for example, that the coefficients $a_{nk}$ obey
the following recursion (deduced from \cite{A1}):
\begin{displaymath}
a_{(n+1)k}=a_{n(k-1)}+2a_{nk}+a_{n(k+1)}.
\end{displaymath}

Expressing the $a_{nk}$'s as in proposition \ref{coefficienti} we
then find:
\begin{displaymath}
[T^k (\beta_{n+1}(x))]_{x=0}=[(T^{k-1}+2T^k +T^{k+1})(\beta_n
(x))]_{x=0},
\end{displaymath}
whence, recalling corollary \ref{powers}:
\begin{equation}\label{ricsomma}
\sum_{h=1}^{n+1}{h-1\choose k-1}f_{(n+1)h}=\sum_{h=1}^{n}\left(
{h-1\choose k-2}+2{h-1\choose k-1}+{h-1\choose k}\right) f_{nh}.
\end{equation}

The sum of binomial coefficients in the r.h.s. of (\ref{ricsomma})
can be easily simplified (using well known properties of the
Pascal matrix) to obtain:
\begin{equation}\label{ricsomma2}
\sum_{h=1}^{n+1}{h-1\choose
k-1}f_{(n+1)h}=\sum_{h=1}^{n}{h+1\choose k}f_{nh}.
\end{equation}

To prove equality (\ref{ricsomma2}) we make use of the structural
properties of the ECO matrix $F$, namely the recursion:
\begin{eqnarray*}
f_{(n+1)h}&=&f_{n(h-1)}+f_{nh}+\ldots +f_{nn}
\\ &=&\sum_{i=h-1}^{n}f_{ni}.
\end{eqnarray*}

Replacing in the l.h.s. of (\ref{ricsomma2}) and interchanging the
order of the summations when necessary, we get:
\begin{eqnarray*}
\sum_{h=1}^{n}{h+1\choose k}f_{nh}&=&\sum_{h=1}^{n+1}{h-1\choose
k-1}\left( \sum_{i=h-1}^{n}f_{ni}\right)
\\ &=&\sum_{i=0}^{n}\left( \sum_{h=1}^{i+1}{h-1\choose k-1}\right)
f_{ni}=\sum_{i=0}^{n}{i+1\choose k}f_{ni},
\end{eqnarray*}
which is an identity ($f_{n0}=0$ by convention). Therefore, we
have formally proved that \emph{switching from the canonical basis
to the Aigner basis translates the ECO matrix $F$ of Catalan
numbers into the (unique) admissible matrix $A$ of Catalan
numbers}.

\subsection{Motzkin numbers}

We can use the same approach to deal with Motzkin numbers.
Consider the succession rule
\begin{displaymath}
\Omega: \left\{ \begin{array}{lll} (1)
\\ (1)\rightsquigarrow (2)
\\ (k)\rightsquigarrow (1)(2)(3)\cdots (k-1)(k+1)
\end{array}
\right. .
\end{displaymath}

$\Omega$ determines the sequence of Motzkin numbers
(\cite{BDLPP}), and its ECO matrix is the following:
\begin{displaymath}
F=\left( \begin{array}{lllllll} 1&0&0&0&0&0&\cdots
\\ 0&1&0&0&0&0&\cdots
\\ 1&0&1&0&0&0&\cdots
\\ 1&2&0&1&0&0&\cdots
\\ 3&2&3&0&1&0&\cdots
\\ 6&7&3&4&0&1&\cdots
\\ \vdots &\vdots &\vdots &\vdots &\vdots &\vdots
&\ddots
\end{array}\right) .
\end{displaymath}

Applying the fundamental change of basis (or, equivalently,
multiplying on the right by the Pascal matrix $P$) we get to the
matrix
\begin{displaymath}
A=\left( \begin{array}{lllllll} 1&0&0&0&0&0&\cdots
\\ 1&1&0&0&0&0&\cdots
\\ 2&2&1&0&0&0&\cdots
\\ 4&5&3&1&0&0&\cdots
\\ 9&12&9&4&1&0&\cdots
\\ 21&30&25&14&5&1&\cdots
\\ \vdots &\vdots &\vdots &\vdots &\vdots &\vdots
&\ddots
\end{array}\right) .
\end{displaymath}

It can be proved that $A$ is the unique admissible matrix
associated with Motzkin numbers. However, it will be an immediate
consequence of the results of section \ref{fsr}.

\subsection{Bell numbers}

The most popular succession rule giving rise to the Bell numbers
is the following:
\begin{displaymath}
\Omega: \left\{ \begin{array}{lll} (1)
\\ (1)\rightsquigarrow (2)
\\ (k)\rightsquigarrow (k)^{k-1}(k+1)
\end{array}
\right. .
\end{displaymath}

Such a rule describes the usual construction of set partitions,
and its ECO matrix is the following:
\begin{displaymath}
F=\left( \begin{array}{llllllll} 1&0&0&0&0&0&0&\cdots
\\ 0&1&0&0&0&0&0&\cdots
\\ 0&1&1&0&0&0&0&\cdots
\\ 0&1&3&1&0&0&0&\cdots
\\ 0&1&7&6&1&0&0&\cdots
\\ 0&1&15&25&10&1&0&\cdots
\\ 0&1&31&89&65&15&1&\cdots
\\ \vdots &\vdots &\vdots &\vdots &\vdots &\vdots &\vdots
&\ddots
\end{array}\right) .
\end{displaymath}

The matrix $F$ is the well-known matrix of the Stirling numbers of
the second kind. Applying the fundamental change of basis leads to
the matrix
\begin{displaymath}
A=\left( \begin{array}{llllll} 1&0&0&0&0&\cdots
\\ 1&1&0&0&0&\cdots
\\ 2&3&1&0&0&\cdots
\\ 5&10&6&1&0&\cdots
\\ 15&37&31&10&1&\cdots
\\ \vdots &\vdots &\vdots &\vdots &\vdots
&\ddots
\end{array}\right) .
\end{displaymath}

It is immediately seen that $A$ is not an admissible matrix.
Nevertheless, the elements of the main diagonal are all equal to
1, so $A$ may be an Aigner matrix. Indeed, it can be shown that
$A$ is the Aigner matrix of type $(\sigma ,\tau)$, where $\sigma
=(k)_{\kinNstar}$ and $\tau =(k-1)_{k\geq 2}$ \cite{A3}. This is
our first example of an Aigner matrix which is not admissible and
is linked to an ECO matrix by the fundamental change of basis.

\subsection{Factorial numbers}

The case of factorial numbers, which is extremely simple from the
point of view of succession rules, turns out to be rather curious
when the fundamental change of basis is applied. Indeed, the
trivial rule
\begin{displaymath}
\Omega: \left\{ \begin{array}{ll} (1)
\\ (k)\rightsquigarrow (k+1)^k
\end{array}
\right.
\end{displaymath}
for the factorial numbers leads to the diagonal ECO matrix
\begin{displaymath}
F=\left( \begin{array}{llllll} 1&0&0&0&0&\cdots
\\ 0&1&0&0&0&\cdots
\\ 0&0&2&0&0&\cdots
\\ 0&0&0&6&0&\cdots
\\ 0&0&0&0&24&\cdots
\\ \vdots &\vdots &\vdots &\vdots &\vdots
&\ddots
\end{array}\right) ,
\end{displaymath}
where, clearly, $f_{nk}=(n-1)!\delta_{nk}$ ($\delta$ is the usual
Kronecker delta). Thus we have $\beta_n (x)=(n-1)!x^n
=\sum_{k=1}^{n}(n-1)!{n-1\choose k-1}p_k (x)$, whence
\begin{displaymath}
A=\left( \begin{array}{llllll} 1&0&0&0&0&\cdots
\\ 1&1&0&0&0&\cdots
\\ 2&4&2&0&0&\cdots
\\ 6&18&18&6&0&\cdots
\\ 24&96&144&96&24&\cdots
\\ \vdots &\vdots &\vdots &\vdots &\vdots
&\ddots
\end{array}\right) .
\end{displaymath}

It is clear that $A$ cannot be an Aigner matrix, since the
elements on the main diagonal are not equal to 1. However, if we
consider the scalar multiplication of the $(n+1)$-st and
$(m+1)$-st rows of $A$, supposing that $n\leq m$, we get:
\begin{eqnarray*}
\sum_{k=1}^{n+1}a_{(n+1)k}a_{(m+1)k}&=&\sum_{k=0}^{n}n!m!{n\choose
k}{m\choose k}
\\ &=^{(*)}&(n+m)!{n+m\choose n}=a_{(n+m+1)1}.
\end{eqnarray*}

(Equality $^{(*)}$ is an application of Vandermonde's
convolution). Thus $A$ possesses a typical property of admissible
matrices, without being neither admissible nor Aigner.

\subsection{Involutions}

Involutions are considered (from a succession rule point of view)
in \cite{FP}. They are generated by the following succession rule:
\begin{displaymath}
\Omega: \left\{ \begin{array}{lll} (1)
\\ (1)\rightsquigarrow (2)
\\ (k)\rightsquigarrow (k-1)^{k-1}(k+1)
\end{array}
\right. ,
\end{displaymath}
giving rise to the ECO matrix:
\begin{displaymath}
F=\left( \begin{array}{llllll} 1&0&0&0&0&\cdots
\\ 0&1&0&0&0&\cdots
\\ 1&0&1&0&0&\cdots
\\ 0&3&0&1&0&\cdots
\\ 3&0&6&0&1&\cdots
\\ \vdots &\vdots &\vdots &\vdots &\vdots
&\ddots
\end{array}\right) .
\end{displaymath}

In this case, the fundamental change of basis leads to the matrix:
\begin{displaymath}
A=\left( \begin{array}{llllll} 1&0&0&0&0&\cdots
\\ 1&1&0&0&0&\cdots
\\ 2&2&1&0&0&\cdots
\\ 4&6&3&1&0&\cdots
\\ 10&16&12&4&1&\cdots
\\ \vdots &\vdots &\vdots &\vdots &\vdots
&\ddots
\end{array}\right) ,
\end{displaymath}
which is the Aigner matrix of type $(\sigma ,\tau )$, for $\sigma
=(1)_{\kinNstar}$ and $\tau =(k-1)_{k\geq 2}$ \cite{A3}. This case
has some analogies with that of Bell numbers (for example, $A$ is
Aigner but not admissible).

\section{Factorial succession rules}\label{fsr}

Referring to \cite{Betal,FP}, we recall the definition of a
factorial succession rule and a factorial rule operator.

A \emph{factorial succession rule} is a rule of the form:
\begin{displaymath}
\Omega: \left\{ \begin{array}{ll} (a)
\\ (k)\rightsquigarrow (r_0 )(r_0 +1)\cdots (r_0 +k-m-1)(k+d_1 )\cdots (k+d_m )
\end{array}
\right. ,
\end{displaymath}
for $k\geq r_0 \geq 1$. A \emph{factorial rule operator} is the
rule operator of a factorial rule. In \cite{Betal} it is shown
that factorial rules have an algebraic generating function. In
\cite{FP} it is stated the following result, concerning the form
of a factorial rule operator.

\begin{prop} (\cite{FP}) A rule operator $L$ is factorial if and only if
$L=p(\mathbf{x},\mathbf{x}^{-1},T)$, where $p(a,b,c)$ is a
polynomial of degree 1 in $c$ having the form:
\begin{displaymath}
p(a,b,c)=u_0 (a)+v_0 (b)+u_1 (a)c,
\end{displaymath}
and $T$ is the factorial derivative operator, as usual.
\end{prop}

In the present section we give a complete answer to the first of
the problems stated in section \ref{bascha} for the class of ECO
matrices arising from factorial rules. In order to accomplish our
result, we need to slightly generalize Aigner's original setting.

Consider an infinite lower triangular matrix
$A=(a_{nk})_{\nkinNstar}$, with $a_{11}=1$, and denote by $L$ the
linear operator associated with its rows, as in (\ref{operiga}).
We say that A is a \emph{generalized Aigner matrix} when there
exist three nonnegative integer sequences $(r_n )_{\ninN}$, $(s_n
)_{\ninN}$, $(t_n )_{\ninN}$ such that, for every $\ninN$:
\begin{equation}\label{ultima}
L(p_n (x))=t_n p_{n-1}(x)+s_n p_n (x)+r_n p_{n+1}(x).
\end{equation}

The following fact follows immediately from the above definition.

\begin{prop}\label{ricor}
If $A$ is a generalized Aigner matrix, then its entries obey the
following recursion:
\begin{equation}\label{genrec}
\left\{ \begin{array}{ll} a_{11}=1
\\ a_{(n+1)k}=r_{k-1}a_{n(k-1)}+s_k a_{nk}+t_{k+1}a_{n(k+1)}
\end{array} \right. .
\end{equation}
\end{prop}

\emph{Proof.} Since $L(\beta_n (x))=\beta_{n+1}(x)$, using
(\ref{ultima}) we have:
\begin{eqnarray*}
\sum_{k=1}^{n+1}a_{(n+1)k}p_k (x)&=&\sum_{k=1}^{n}a_{nk}L(p_k (x))
\\ &=&\sum_{k=1}^{n}a_{nk}\left( t_k p_{k-1}(x)+s_k p_k
(x)+r_k p_{k+1}(x)\right) \\ &=&\sum_{k=1}^{n+1}\left(
r_{k-1}a_{n(k-1)}+s_k a_{nk}+t_{k+1}a_{n(k+1)}\right) p_k (x),
\end{eqnarray*}
whence the thesis follows.\cvd

In particular, it is easily seen that a generalized Aigner matrix
is not forced to have the elements on the main diagonal equal to
1.

\bigskip

\textbf{Question}: \emph{does the converse of proposition
\ref{ricor} hold?}

\bigskip

Now we can focus on the case of factorial rules. Suppose that $L$
is a factorial rule operator of the form
\begin{equation}\label{factop}
L=a(\mathbf{x})+b(\mathbf{x}^{-1})+c(\mathbf{x})T,
\end{equation}
where $a(x)=\sum_k a_k x^k$, $b(x)=\sum_{k\geq 1}b_k x^{_k}$,
$c(x)=\sum_k c_k x^k$ are fixed polynomials. We can immediately
find a sufficient condition for $L$ to induce a generalized Aigner
matrix.

\begin{prop} If $b(x^{-1})=0$, deg $a(x)\leq 1$ and deg $c(x)\leq 2$,
then $A$ is a generalized Aigner matrix.
\end{prop}

\emph{Proof.} We have to show that $L$ acts as in (\ref{ultima}).
Indeed, a simple computation shows that
\begin{eqnarray*}
L(p_n (x))&=&((a_0 +a_1 x)+(c_0 +c_1 x+c_2 x^2)T)(p_n (x))
\\ &=&(c_0 +c_1 +c_2 )p_{n-1}(x)+(a_0 +a_1 +c_1 +2c_2 )p_n (x)+(a_1
+c_2 )p_{n+1}(x),
\end{eqnarray*}
which is enough to conclude thanks to the previous
proposition.\cvd

\bigskip

\emph{Examples.} \begin{itemize} \item[i)] The succession rules
for Catalan and Motzkin numbers described above are associated
with rule operators for which, respectively,
$a(x)=b(x)=0,c(x)=x^2$ and $a(x)=x-1,b(x)=0,c(x)=x$. \item[ii)]
Consider the following succession rule, inducing Schr\"oder
numbers:
\begin{equation}\label{schroder}
\left\{ \begin{array}{lll} (1)
\\ (1)\rightsquigarrow (2)
\\ (k)\rightsquigarrow (3)(4)(5)\cdots (k-1)(k)(k+1)^2
\end{array}
\right. .
\end{equation}

In this case, the rule operator $L$ has the form:
\begin{displaymath}
L=\mathbf{x}-\mathbf{x}^2 +\mathbf{x}^3 T,
\end{displaymath}
so it does not satisfy the hypotheses of the above proposition.
Nevertheless, this rule is related to a generalized Aigner matrix,
as we will see in the next pages.
\end{itemize}

\begin{teor} A factorial rule operator $L$ as in (\ref{factop}) is
associated with a generalized Aigner matrix if and only if the
following conditions hold:
\begin{itemize}
\item[i)] $\sum_k {k\choose h}a_k +\sum_k {k\choose h+1}c_k
=0,\qquad \forall h\geq 2$;
\item[ii)] $\sum_{k\geq
1}(-1)^{n-k-h}{n-1-h\choose k-1}b_k =0,\qquad \forall h<n-1$.
\end{itemize}
\end{teor}

\emph{Proof.} The rule operator $L$ acts on the Aigner basis as
follows:
\begin{eqnarray*}
L(p_n (x))&=&\sum_k a_k x^k p_n (x)+\sum_{k\geq 1}b_k x^{-k}p_n
(x)+\sum_k c_k x^k p_{n-1}(x)
\\ &=&\sum_k a_k \left( \sum_{h=0}^{k}{k\choose h}p_{n+h}(x)\right)
\\ &+&\sum_{k\geq 1}b_k \left( \sum_{h=0}^{n-k}(-1)^{n-k-h}{n-1-h\choose
k-1}p_h (x)\right)
\\ &+&\sum_k c_k \left( \sum_{h=0}^{k}{k\choose h}p_{n-1+h} (x)\right) .
\end{eqnarray*}

Now, $L$ must satisfy condition (\ref{ultima}), which means that,
in the above expansion, all the coefficients of the polynomials
$p_k (x)$, for $k\notin \{ n-1,n,n+1\}$, must be zero. This
translates into conditions $i)$ and $ii)$ above, so the proof is
complete.\cvd

\bigskip

\emph{Example.} The previous example related to Schr\"oder numbers
can now be reconsidered. It is clear that condition $ii)$ is
trivially verified, whereas the only interesting case of condition
$i)$ occurs when $h=2$, and we have:
\begin{displaymath}
{2\choose 2}(-1)+{3\choose 3}1=-1+1=0.
\end{displaymath}

So rule (\ref{schroder}) is associated with a \emph{generalized}
Aigner matrix $A$. By an explicit computation for the operator
$L$, we find for $A$ the following expression:
\begin{displaymath}
A=\left( \begin{array}{llllll} 1&0&0&0&0&\cdots
\\ 1&1&0&0&0&\cdots
\\ 2&4&2&0&0&\cdots
\\ 6&16&14&4&0&\cdots
\\ 22&68&78&40&8&\cdots
\\ \vdots &\vdots &\vdots &\vdots &\vdots
&\ddots
\end{array}\right) ,
\end{displaymath}
where the entries obey the following recursion:
\begin{displaymath}
\left\{ \begin{array}{lll} a_{(n+1)1}=a_{n1}+a_{n2},
\\ a_{(n+1)2}=a_{n1}+3a_{n2}+a_{n3},
\\ a_{(n+1)k}=2a_{n(k-1)}+3a_{nk}+a_{n(k+1)},\qquad \textnormal{for $k>2$}
\end{array}
\right. .
\end{displaymath}

At this stage, it is worth noting that in \cite{A3} an Aigner
matrix for Schr\"oder numbers is taken into consideration,
precisely:
\begin{displaymath}
\tilde{A}=\left( \begin{array}{llllll} 1&0&0&0&0&\cdots
\\ 2&1&0&0&0&\cdots
\\ 6&5&1&0&0&\cdots
\\ 22&23&8&1&0&\cdots
\\ 90&107&49&11&1&\cdots
\\ \vdots &\vdots &\vdots &\vdots &\vdots
&\ddots
\end{array}\right) .
\end{displaymath}

If we apply the inverse of the fundamental change of basis (that
is, we switch from $p_n (x)$ to $x^n$), we find the following
matrix:
\begin{displaymath}
\tilde{F}=\left( \begin{array}{llllll} 1&0&0&0&0&\cdots
\\ 1&1&0&0&0&\cdots
\\ 2&3&1&0&0&\cdots
\\ 6&5&10&1&0&\cdots
\\ 22&38&22&7&1&\cdots
\\ \vdots &\vdots &\vdots &\vdots &\vdots
&\ddots
\end{array}\right) .
\end{displaymath}

Strictly speaking, $\tilde{F}$ is not an ECO matrix, as we can
immediately notice. However, if we suppose that, in $\tilde{F}$,
column $k$ represents the distribution of the label $(2k)$, we
manage to find an ECO-interpretation. From a linear operator point
of view, we can consider the operator $\tilde{L}$ associated with
the rows of $\tilde{F}$: if we replace the variable $x$ with
$x^2$, we in fact obtain a rule operator $L$, which corresponds to
a well-known ECO-interpretation of Schr\"oder numbers
\cite{BDLPP2}. The succession rule related to $L$ is the
following:
\begin{displaymath}
\left\{ \begin{array}{ll} (2)
\\ (2k)\rightsquigarrow (2)(4)^2 (6)^2 \cdots (2k)^2 (2k+2)
\end{array}
\right. .
\end{displaymath}

\section{Differential succession rules}

We call \emph{differential succession rule} each rule such that in
the production of every label $(k)$ at least one label greater or
equal than $(k-1)$ has an exponent linearly depending on $(k)$.
Equivalently, a \emph{differential rule operator} is a rule
operator which can be expressed in the form
\begin{displaymath}
L=p(\mathbf{x},\mathbf{x}^{-1},D)=a(\mathbf{x})+b(\mathbf{x}^{-1})+c(\mathbf{x})D.
\end{displaymath}

Observe that, for reasons of consistency (in a succession rule a
node labelled $(k)$ must produce exactly $k$ sons), in the above
formula we necessarily have $c(x)=x^t$, for some $t\in
\mathbf{N}$, so that a differential rule operator has the
following, general expansion:
\begin{equation}\label{dro}
L=p(\mathbf{x},\mathbf{x}^{-1},D)=a(\mathbf{x})+b(\mathbf{x}^{-1})+\mathbf{x}^t
D.
\end{equation}

The rules given in section \ref{esempi} for Bell numbers,
factorials and the number of involutions of a set are examples of
differential succession rules, whose associated differential rule
operators are, respectively,
$\mathbf{x}-1+\mathbf{x}D,\mathbf{x}^2 D$ and
$\mathbf{x}-\mathbf{x}^{-1}+D$. In \cite{FP} it is shown that the
sequences determined by a differential succession rule possess a
transcendental (ordinary) generating function. Using arguments
which are essentially analogue to those employed for factorial
succession rules, we can prove the following results, concerning
the relationship with Aigner's theory.

\begin{teor}\label{cardro} A differential rule operator $L$ as in (\ref{dro}) is
associated with a generalized Aigner matrix if and only if the
following conditions hold:
\begin{enumerate}
\item if $t\neq 0$, the conditions are:
\begin{itemize}
\item[i)] $\sum_k {k\choose h}a_k +n{t\choose h+1}-{t-1\choose
h+1}=0,\qquad \forall h\geq 2$; \item[ii)] $\sum_{k\geq
1}(-1)^{n-k-h}{n-1-h\choose k-1}b_k=0,\qquad \forall h<n-1$;
\end{itemize}
\item if $t=0$, the conditions are:
\begin{itemize}
\item[i)] $\sum_k {k\choose h}a_k =0,\qquad \forall h\geq 2$;
\item[ii)] $\sum_{k\geq 1}(-1)^{n-k-h}{n-1-h\choose k-1}b_k
-(-1)^{n-h}=0,\qquad \forall h<n-1$;
\end{itemize}
\end{enumerate}
\end{teor}

\begin{cor} If $b(x^{-1})=0$, deg $a(x)\leq 1$ and $t\in \{
1,2\}$, then $L$ is associated with a generalized Aigner matrix.
\end{cor}

\emph{Examples.} The cases of Bell numbers and factorial numbers
can be easily tackled using the corollary (see above for the rule
operators involved). As far as involutions are concerned, we have
to consider the rule operator $\mathbf{x}-\mathbf{x}^{-1}+D$.
Applying theorem \ref{cardro} in the case $t=0$, we have that
condition \emph{i)} is trivially satisfied, whereas condition
\emph{ii)} becomes:
\begin{displaymath}
(-1)^{n-1-h}{n-1-h\choose
0}(-1)-(-1)^{n-h}=(-1)^{n-h}-(-1)^{n-h}=0,
\end{displaymath}
which is enough to conclude.

\section{Conclusions and further work}

The present work is intended to be only the first step towards a
more detailed investigation of the relationship between the ECO
method and the theory of Catalan-like numbers.

One of the first things to be done in the next future is to
provide a more complete gallery of examples to illustrate the
soundness of our approach. In particular, it would be nice to find
new applications of the ECO method starting from known ones of
Aigner's theory, and vice versa, which is what we hope to do in a
forthcoming publication.

Another line of research is provided by the last section of
\cite{A3}, where Aigner introduces the basics of what he calls
``ballot enumeration". The analogies with the techniques of the
ECO method are evident and, in fact, the combinatorial model
proposed by Aigner is a particular instance of an ECO
construction.

\end{document}